\documentclass[letterpaper, 10 pt, conference]{ieeeconf}  

\usepackage{cite}
\usepackage{graphicx}
\usepackage{amsmath}
\usepackage{amssymb}
\allowdisplaybreaks[4]
\usepackage{mathrsfs}
\usepackage{verbatim}
\usepackage{amsfonts}
\usepackage{float}
\usepackage{subfigure}
\usepackage{caption}
\usepackage{hhline}
\usepackage{delarray}
\usepackage{epstopdf}
\usepackage{setspace}
\usepackage{enumerate}
\usepackage{lineno}
\usepackage[colorlinks=true, allcolors=blue]{hyperref}
\usepackage{color}

\IEEEoverridecommandlockouts                  

\title{\LARGE \bf
Region of Synchronization Estimation for Complex Networks \\
via SOS Programming
}

\author{Shuyuan~Zhang$^{1}$, Rapha\"{e}l M. Jungers$^{1}$, and Lei Wang$^{2}$
\thanks{RJ is a FNRS honorary Research Associate. This project has received funding from the \textit{European Research Council (ERC) under the European Union's Horizon 2020 research and innovation programme} under grant agreement No 864017 - L2C, from the Horizon Europe programme under grant agreement No101177842 - Unimaas, and from the ARC (French Community of Belgium) - project name: SIDDARTA.}
\thanks{$^{1}$S. Zhang and R. M. Jungers are with the ICTEAM Institute, UCLouvain, 4 Avenue Georges Lema\^itre, 1348 Louvain-la-Neuve, Belgium. {\tt\small shuyuan.zhang,raphael.jungers@uclouvain.be}} 
\thanks{$^{2}$L. Wang is with the School of Automation Science and Electrical Engineering, Beihang University, Beijing 100191, China. {\tt\small lwang@buaa.edu.cn}}
}

\begin{document}

\maketitle
\thispagestyle{empty}
\pagestyle{empty}

\begin{abstract}

In this article, we explore the problem of the region of synchronization (ROS) for complex networks with nonlinear dynamics. Given a pair of state- and target- sets, our goal is to estimate the ROS such that the trajectories originating within it reach the target set (i.e., synchronization manifold), without leaving the state set before the first hitting time. In order to do so, an exponential guidance-barrier function is proposed to construct the ROS along the synchronization manifold, and the corresponding sufficient conditions for estimating the ROS are developed. The resulting conditions lead to a sum-of-squares programming problem, thereby affording a polynomial-time solvability.
Furthermore, when the synchronization manifold reduces to an equilibrium point, our method not only estimates a larger ROS compared to existing results but also allows the ROS to take more general shapes.
Finally, we present two numerical examples to demonstrate the effectiveness of the theoretical results.

\end{abstract}

\vspace{0.2em}
\section{Introduction}
\label{1}

Complex networks are composed of nodes and edges, where the nodes represent individual agents and the edges capture the interactions between them.
Due to their ability to model real-world interconnected systems, complex networks have served as a ubiquitous paradigm with a broad range of applications across fields such as physics, biology, economics, control science and engineering \cite{berner2021desynchronization,strogatz2001exploring,arenas2008synchronization,d2016complex,modares2018optimal,sun2023mean}.
Synchronization in complex networks, recognized as a fundamental collective phenomenon in populations of interconnected nodes, has received significant attention from the scholarly community. A major question among researchers is how to achieve synchronization in complex dynamical networks characterized by nonlinear node dynamics and interactive couplings.
A prominent work is the master stability function method proposed in \cite{pecora1998master}, which investigated the stability of synchronous states through the use of Lyapunov exponents.
The subsequent literature has focused on establishing synchronization criteria related to coupling strength for complex networks, addressing both linear and nonlinear dynamics and coupling structures \cite{Zhang2:2021,ma2010necessary,delellis2011quad,zhang2022polynomial,zhang2024automatic,panteley2017synchronization,wang2021decomposition}.

In addition to establishing synchronization criteria, another critical aspect is identifying the initial conditions from which the trajectories of complex networks converge to the synchronous states as time progresses. The set of initial conditions from which synchronization is achieved is referred to as the region of synchronization (ROS).
The concept of ROS is analogous to the well-known region of attraction (ROA); however, they differ significantly. The ROA pertains solely to the convergence of the system's trajectories to an equilibrium point, whereas the ROS exists along a synchronization manifold.
Regardless of which one is considered, the research objective remains the same: to estimate them as accurately as possible, since the exact computation of either is often challenging.
In \cite{chesi2011domain}, Chesi presented a unified framework for estimating the ROA of nonlinear dynamical systems, including polynomial, uncertain polynomial, and non-polynomial systems, using the emerging sum-of-squares (SOS) programming approach, which involves solving convex optimization problems via semi-definite programming.
In \cite{iannelli2019region}, an SOS-based optimization algorithm was proposed to enlarge the inner estimates of ROA with integral quadratic constraints.
In \cite{korda2013inner} and \cite{henrion2014convex}, Korda \emph{et al.} developed an infinite-dimensional convex characterization of ROA for polynomial systems,
providing inner and outer approximations of the ROA, respectively.
Recently, Xue \emph{et al.} \cite{xue2023RA-Verification} discussed the under-approximation of the reach-avoid set based on SOS decomposition techniques. When the safe set is defined as the state set and the target set is an equilibrium point, the reach-avoid set corresponds to the ROA.

However, there is little research on ROS for complex networks, particularly regarding computational approaches to its estimation.
As is known, the shape of the ROS is generally quite complicated due to the node dynamics and the coupling relationships within the networks.
For linear networks, the ROS in \cite{zhang2011optimal} and \cite{li2010consensus} is defined on the complex plane formed by the product of the coupling strength and the eigenvalues of the Laplacian matrix.
In \cite{guaracy2023insights}, Guaracy \emph{et al.} further explored the selection of the coupling and feedback gains for estimating the ROS, drawing on robust control and small-gain theorem perspectives. 
For nonlinear networks, the authors in \cite{liu2007analyzing, duan2009disconnected, gequn2014control} discussed the sufficient conditions for achieving either an unbounded or a bounded synchronization region by linearizing the node dynamics around an equilibrium point. 
In \cite{zhu2019estimating} and \cite{zhu2023estimating}, Zhu \emph{et al.} estimated the ROS of an equilibrium point in a nonlinear dynamical network by computing quadratic Lyapunov functions, and analyzed the network stability.
However, the above studies primarily focused on the problem of ROS around an equilibrium point. Moreover, the estimated ROS is often conservative and tends to have a simple shape, such as a circle or ellipse \cite{zhu2019estimating}. 
Thus, this article aims to present a method for estimating the ROS of nonlinear networks with nonlinear couplings along the synchronization manifold, which may represent an equilibrium point, a periodic trajectory, or even a chaotic trajectory, while capturing a more general shape of the ROS.

In this article, we present an exponential guidance-barrier function (EGBF) to construct the ROS along the synchronization manifold for nonlinear networks. The trajectories starting from the ROS are ensured to reach the synchronization manifold within a finite time.
For polynomial networks, the constraints resulting from the EGBF fall within a framework of convex SOS programming, which can be efficiently solved by semi-definite programming in polynomial time.
In particular, when the synchronization manifold reduces to an equilibrium point, our method not only estimates a larger ROS than \cite{zhu2019estimating} but also allows the ROS to take more general shapes.
Finally, two numerical examples are given to demonstrate the effectiveness of our method.

The rest of this article is organized as follows. Section \ref{2} formulates the problem of interest, and some mathematical preliminaries. In Section \ref{3}, the ROS are constructed by EGBF and estimated by SOS programming. After demonstrating our approach through two numerical examples in Section \ref{4}, we wrap up this article in Section \ref{5}.

Notations: $\mathcal{N}$ is the set of integer numbers $\{1,\ldots,N\}$,
$\mathbb{R}$ is the set of real numbers, $\mathbb{R}^{n}$ is the set of real vectors of $n$-dimension, $\mathbb{R}^{n \times n}$ is the set of $n\times n$ dimensional real matrices, and $\mathbb{R}[\cdot]$ is the set of real polynomials.
The closure of a set $\mathcal{M}$ is denoted by $\overline{\mathcal{M}}$, and the boundary by $\partial\mathcal{M}$.
The Euclidean norm of vector $\delta$ is denoted by $||\delta||$, and the transpose of $\delta$ is denoted by $\delta^\top$.
For matrices $\Phi$ and $\Psi$, the Kronecker product is represented by $\Phi\otimes \Psi$.
$\sum[x]$ denotes the set of multivariate SOS polynomials, i.e.,
\[\sum[x] = \left \{\rho(x) \in \mathbb{R}[x]  \middle |\;
\begin{aligned}
 &\rho(x)=\sum_{i=1}^k p_i^2(x), \forall x \in \mathbb{R}^{n} \\
 &p_i(x)\in \mathbb{R}[x], i=1,\ldots,k
\end{aligned}
\right\}.
\]

\vspace{0.5em}
\section{Preliminaries}
\label{2}

\subsection{Region of Attraction}
\label{2-1}

In this subsection, we review the concept of ROA for nonlinear systems.

Consider an autonomous system in the form,
\begin{align}
\dot{x}=f(x),  \label{autonomous}
\end{align}
where $x \in \mathbb{R}^{n}$ is the state vector, and $f(\cdot)$ is locally Lipschitz continuous.

Let $\varphi_{x_{0}}(\cdot): [0, T_{x_{0}}) \rightarrow \mathbb{R}^{n}$ denote the trajectory of system \eqref{autonomous} originating from $x_{0} \in \mathbb{R}^{n}$ over the time domain $[0, T_{x_{0}})$. Then, we have 
\begin{align*}
\varphi_{x_{0}}(t):=x(t), \ \forall t \in [0,T_{x_{0}})
\end{align*}
with $\varphi_{x_{0}}(0)=x_{0}$, and $T_{x_{0}}$ is a positive constant.

Given a bounded open state set $\mathcal{X}$ and a target set $\mathcal{X}_{T}$ in $\mathcal{X}$, where 
\begin{enumerate}[(1)]
     \item $\mathcal{X}=\{x \in \mathbb{R}^{n}| h(x)>0\}$,
     \item $\overline{\mathcal{X}}=\{x \in \mathbb{R}^{n}| h(x)\geq 0\}$,
     \item $\partial\mathcal{X}=\{x \in \mathbb{R}^{n}| h(x)=0\}$,
     \item $\mathcal{X}_{T}=\{x \in \mathbb{R}^{n}| ||x||< \varepsilon \}$,
\end{enumerate}
where $\varepsilon$ is a positive constant, and $h(\cdot): \mathbb{R}^{n} \rightarrow \mathbb{R}$ is a continuous function.

The ROA property with respect to them is elaborated in the following definition.

{\bf Definition 1} \cite{korda2013inner}.
Given a pair of state- and target- sets $\mathcal{X}$ and $\mathcal{X}_{T}$, if there exists a $T>0$ such that for any trajectory $x(t)$ of system \eqref{autonomous} starting from $x_{0}$, there holds
\begin{align}
\varphi_{x_{0}}(T) \in \mathcal{X}_{T} \bigwedge \varphi_{x_{0}}(t) \in \mathcal{X} \ \text{for} \ t \in [0,T], \label{RA}
\end{align}
then we say that $\mathcal{X}_{0}=\{x_{0} \in \mathbb{R}^{n} | \ \eqref{RA} \ \text{holds}\}$ is the ROA of system \eqref{autonomous}.

Given that the origin is a stable equilibrium point within the target set, we define the ROA $\mathcal{X}_{0}$. This is the set of all initial conditions for which there exists an admissible trajectory, i.e., the set of all initial conditions that can be steered to the target set $\mathcal{X}_{T}$ eventually while staying inside the state set $\mathcal{X}$ until the first time the target is reached.

\subsection{Problem Formulation}
\label{2-2}
In this subsection, we present the problem of ROS along the synchronization manifold for complex networks.

Consider a complex network consisting of $N$ nonlinearly coupled nodes, described by
\begin{align} \label{network}
\dot{\mathbf{x}}_{i}=f(\mathbf{x}_{i})-c\sum_{j=1}^{N} L_{ij} g(\mathbf{x}_{j}), ~i \in \mathcal{N},
\end{align}
where $\mathbf{x}_i=[x_{i1},x_{i2},\ldots,x_{in}]^{\top}\in\mathbb{R}^n$ is the state vector of the $i$-th node, $f(\mathbf{x}_{i})=[f_{1}(x_{i}),f_{2}(x_{i}),\ldots,f_{n}(x_{i})] \in \mathbb{R}^n \rightarrow \mathbb{R}^n$ is the smooth vector field, $c>0$ is the coupling strength, $g(\mathbf{x}_{j})=[g_{1}(x_{j}),g_{2}(x_{j}),\ldots,g_{n}(x_{j})] \in \mathbb{R}^n \rightarrow \mathbb{R}^n$ is the nonlinear coupling function.
The communication topology among the $N$ nodes is represented by a digraph $\mathcal{G}$. The corresponding Laplacian matrix $L=(L_{ij})\in\mathbb{R}^{N\times N}$ is defined as follows. If there is a directed connection from node $j$ to node $i$~($i\neq j$), then $L_{ij}<0$; otherwise, $L_{ij}=0$; and the diagonal entries of the matrix $L$ are defined by $L_{ii}=-\sum^N_{j=1,j\neq i}L_{ij}$, $i \in \mathcal{N}$.

Let $\mathbf{X}=[\mathbf{x}_{1}^{\top},\mathbf{x}_{2}^{\top},\ldots,
\mathbf{x}_{N}^{\top}]^{\top} \in \mathbb{R}^{nN}$ be the aggregate states of all $N$ nodes of network \eqref{network}. The solution trajectory $\mathbf{X}(t;\mathbf{X}_{0})$ of \eqref{network} for any initial state $\mathbf{X}(0)=\mathbf{X}_{0} \in \mathbb{R}^{nN}$ is defined over the time interval $[0,+\infty)$, and may be simply denoted as $\mathbf{X}(t)$ or even $\mathbf{X}$.

{\bf Definition 2} \cite{wang2021decomposition}. The set $\mathcal{S}=\{\mathbf{X}=[\mathbf{x}_1^{\top},\mathbf{x}_2^{\top},\ldots,\mathbf{x}_N^{\top}]^{\top}$ $ \in \mathbb{R}^{nN}|~\mathbf{x}_i=\mathbf{x}_j,~i,j \in \mathcal{N}\}$ is called the synchronization manifold of network \eqref{network}.

In the synchronization process, the trajectories of network \eqref{network} may converge to an equilibrium point, i.e., $\mathbf{x}_i  \rightarrow \mathbf{x}_j  \rightarrow \gamma$ (where $\gamma$ is a constant vector), or even they may converge to a closed orbit, i.e., $\mathbf{x}_i  \rightarrow \mathbf{x}_j  \rightarrow \mathbf{x}^{*}(t)$ (where $\mathbf{x}^{*}(t)$ is a time-varying trajectory). 
In either case, we can define an error vector as $\triangle \mathbf{x}_{ij}=\mathbf{x}_{i}-\mathbf{x}_{j}$ with $i\neq j$ for $i,j \in \mathcal{N}$. The error vector naturally captures relative differences between nodes, thus simplifying the synchronization analysis \cite{panteley2017synchronization}.

Then, the ROS problem along the synchronization manifold $\mathcal{S}$, with respect to the state set $\mathcal{X}$ and the target set $\mathcal{X}_{T}$, can be formulated as follows.

\textbf{Problem 1}. Given a pair of state- and target- sets $\mathcal{X}$ and $\mathcal{X}_{T}$, if there exists a $T > 0$ such that for any error trajectory $\triangle \mathbf{x}_{ij}(t)$ of  network \eqref{network} starting from $\triangle \mathbf{x}_{ij}^{0}$ satisfies
\begin{align} \label{RA-network}
\varphi_{\triangle \mathbf{x}_{ij}^{0}}(T) \in \mathcal{X}_{T} \bigwedge \varphi_{\triangle \mathbf{x}_{ij}^{0}}(t) \in \mathcal{X} \ \text{for} \ t \in [0,T],
\end{align}
where $i \neq j$ and $i,j \in \mathcal{N}$, then we say that $\mathcal{X}_{0}=\{\triangle \mathbf{x}_{ij}^{0} \in \mathbb{R}^{n} | \ \eqref{RA-network} \ \text{holds}\}$ is the ROS of network \eqref{network}.

{\bf Remark 1}. The ROS denoted as $\mathcal{X}_{0}$ is the set of initial error conditions that will lead to network \eqref{network} reaching the synchronization manifold $\mathcal{S}$ in a finite time.
If $\mathcal{S}$ is a stable synchronization manifold, then we say that trajectories starting from $\mathcal{X}_{0}$ will achieve synchronization.
In particular, if the synchronization manifold $\mathcal{S}$ is merely an equilibrium point,
Eq. \eqref{RA-network} simplifies to
\begin{align} \label{RA-network2}
\varphi_{\mathbf{X}_{0}}(T) \in \mathcal{X}_{T} \bigwedge \varphi_{\mathbf{X}_{0}}(t) \in \mathcal{X} \ \text{for} \ t \in [0,T],
\end{align}
and the ROS is given by $\mathcal{X}_{0}=\{\mathbf{X}_{0} \in \mathbb{R}^{nN} | \ \eqref{RA-network2} \ \text{holds}\}$, which is equivalent to the form in (2.2) of \cite{zhu2019estimating}.
Besides, when network \eqref{network} consists of only one node, the expression of ROS in \eqref{RA-network2} reduces to the expression of ROA in \eqref{RA}.

The objective of this article is to estimate the ROS by finding a largest possible open subset of $\mathcal{X}_{0}$. The proposed method is detailed in the following section.

\vspace{0.5em}
\section{Region of Synchronization Estimation}
\label{3}
In this section, based on the SOS decomposition for multivariate polynomials, we will estimate the ROS along the synchronization manifold and the ROS around the equilibrium point, respectively.

\subsection{ROS along Synchronization Manifold}
\label{3-1}
In this subsection, we introduce the notion of EGBF
to construct the ROS for network \eqref{network}, and solve the corresponding constraints by SOS programming.

{\bf Definition 3}. A continuously differentiable function $V(\cdot):  \overline{\mathcal{X}} \rightarrow \mathbb{R}$ is called EGBF, if there exists $\lambda > 0$ such that for any $\mathbf{x}_{i}, \mathbf{x}_{j} \in \mathbb{R}^{n}$ with $i,j \in \mathcal{N}$, there holds
\begin{align}
\mathcal{L}_{V}(\triangle \mathbf{x}_{ij},\mathbf{x}_{i},\mathbf{x}_{j})-\lambda V(\triangle &\mathbf{x}_{ij}) \geq 0, 
\nonumber\\
&\forall \triangle \mathbf{x}_{ij} \in \overline{\mathcal{X} \setminus \mathcal{X}_{T}}, \label{EGBF1} \\
V(\triangle \mathbf{x}_{ij}) \leq 0, \ \forall \triangle \mathbf{x}_{ij} & \in   \partial \mathcal{X}, \label{EGBF2}
\end{align}
where $\mathcal{L}_{V}(\triangle \mathbf{x}_{ij},\mathbf{x}_{i},\mathbf{x}_{j})=\nabla V(\triangle \mathbf{x}_{ij})(\dot{\mathbf{x}}_{i}-\dot{\mathbf{x}}_{j})$.

Inspired by \cite{xue2023RA-Verification}, we extend the concept of EGBF to complex networks.
Next, the EGBF will be used to construct the ROS for network \eqref{network}, as formulated in Theorem 1.

{\bf Theorem 1}. If $V(\cdot)$ is an EGBF, then $\mathcal{R}=\{\triangle \mathbf{x}_{ij} \in  \mathcal{X} \mid V(\triangle \mathbf{x}_{ij})>0 \}$ ($i \neq j$, $i,j \in \mathcal{N}$) is a subset of ROS along synchronization manifold for network \eqref{network}.

{\bf Proof}. For $\forall \triangle \mathbf{x}_{ij}^{0} \in \mathcal{R}$, we integrate both sides of \eqref{EGBF1}, yielding
\begin{align}
V(\varphi_{\triangle \mathbf{x}_{ij}^{0}}(t)) \geq e^{\lambda t} V(\triangle \mathbf{x}_{ij}^{0}) > 0 \label{TH1-1}
\end{align}
for $t \in [0,T]$, and $V(\varphi_{\triangle \mathbf{x}_{ij}^{0}}(t)) > 0$
indicates that the trajectories starting from $\mathcal{R}$ cannot touch the boundary
$\partial\mathcal{R}$. The boundary condition in \eqref{EGBF2} ensures this as well, even if $\partial\mathcal{R} \cap \partial \mathcal{X} \neq \emptyset$.
Thus, the trajectories will not leave $\mathcal{R}$ before arriving at $\mathcal{X}_{T}$.
Since $V(\cdot)$ is bounded over $\overline{\mathcal{X}}$, the trajectories starting from $\mathcal{R}$ cannot remain indefinitely within $\mathcal{R}$ but outside $\mathcal{X}_{T}$, and will eventually enter $\mathcal{X}_{T}$.
\hfill{$\blacksquare$}

{\bf Remark 2}. According to \eqref{TH1-1}, the parameter $\lambda > 0$ serves as an exponential convergence rate, ensuring the trajectories originating from $\mathcal{R}$ converge exponentially towards $\mathcal{X}_{T}$.
Although a larger $\lambda$ has a faster convergence rate, a smaller $\lambda$ leads to the weaker constraint \eqref{EGBF1}. We advise against using a very small $\lambda$ to avoid potential numerical errors.

Different from the existing literature \cite{liu2007analyzing,duan2009disconnected,gequn2014control},
Theorem 1 generates the ROS $\mathcal{R}$ along the synchronization manifold $\mathcal{S}$, rather than merely around an equilibrium point.
Next, we will present our approach for computing this as-yet-undetermined function $V(\cdot)$ to obtain $\mathcal{R}$. 
At present, there are few effective computational techniques for computing $V(\cdot)$ under multivariate non-negativity constraints. 
In fact, even verifying the non-negativity of a fourth-degree polynomial is an NP-hard problem \cite{papachristodoulou2005analysis}. 
Moreover, the non-negativity constraints in \eqref{EGBF1} and \eqref{EGBF2} are inherently complex.
In what follows, we will focus on determining a polynomial form for $V(\cdot)$ of the polynomial network \eqref{network}, with respect to a pair of semi-algebraic state- and target- sets.

To this end, we encode the non-negativity constraints \eqref{EGBF1} and \eqref{EGBF2} into SOS constraints by using the SOS decomposition technique, leading to a semi-definite program as follows:
\begin{equation} \label{SOS1}
\begin{split}
&\max_{\boldsymbol{c},\lambda} \ \boldsymbol{c}^{\top}  \boldsymbol{\rho} \\
&\ \text{s.t.}  \ \mathcal{L}_{V}(\triangle \mathbf{x}_{ij},\mathbf{x}_{i},\mathbf{x}_{j})-\lambda V(\triangle \mathbf{x}_{ij})
- p_{1}(\triangle \mathbf{x}_{ij})h(\triangle \mathbf{x}_{ij}) \\
&\quad \quad \quad \quad \ \ + p_{2}(\triangle \mathbf{x}_{ij})l(\triangle \mathbf{x}_{ij}) \in \sum[\triangle \mathbf{x}_{ij},\mathbf{x}_{i},\mathbf{x}_{j}],    \\
& \quad \  -V(\triangle \mathbf{x}_{ij})+q(\triangle \mathbf{x}_{ij})h(\triangle \mathbf{x}_{ij}) \in \sum[\triangle \mathbf{x}_{ij}], \\
& \quad \ \ p_{1}(\triangle \mathbf{x}_{ij}) \in \sum[\triangle \mathbf{x}_{ij}],
\ p_{2}(\triangle \mathbf{x}_{ij}) \in \sum[\triangle \mathbf{x}_{ij}], \\
& \quad \ \ q(\triangle \mathbf{x}_{ij}) \in \mathbb{R}[\triangle \mathbf{x}_{ij}],
\end{split}
\end{equation}
where $\boldsymbol{c}^{\top} \boldsymbol{\rho} = \int_{\overline{\mathcal{X}}} V(\triangle \mathbf{x}_{ij})d(\triangle \mathbf{x}_{ij})$, $\boldsymbol{\rho}$ is a constant vector computed by integrating the monomials in $V(\triangle \mathbf{x}_{ij}) \in \mathbb{R}[\triangle \mathbf{x}_{ij}]$ over $\overline{\mathcal{X}}$, and $\boldsymbol{c}$ is a vector composed of unknown coefficients in $V(\triangle \mathbf{x}_{ij})$; and $l(\triangle \mathbf{x}_{ij})=\varepsilon - ||\triangle \mathbf{x}_{ij}||$.

{\bf Remark 3}. In order to obtain a less conservative feasible solution for $V(\cdot)$, we incorporate an objective function in program \eqref{SOS1}, to maximize the integral of $V(\triangle \mathbf{x}_{ij})$ over $\overline{\mathcal{X}}$.
However, it is worthy noting that the volume of $\mathcal{R}$ is not equivalent to the introduced objective function $\boldsymbol{c}^{\top}  \boldsymbol{\rho}$. Since an SOS condition cannot be employed to describe $\mathcal{R}$ to be solved, we instead use $\boldsymbol{c}^{\top}  \boldsymbol{\rho}$ to obtain a numerical approximation for $V(\cdot)$ \cite{xue2023RA-Verification, zhao2023inner}.
Besides, the computation of the constant vector $\boldsymbol{\rho}$ depends on the shape of the given set $\overline{\mathcal{X}}$. If the set $\overline{\mathcal{X}}$ has a simple shape such as circle or ellipse, we can analytically compute the integral based on its boundary conditions \cite{zhao2023inner}. 
In the case where $\overline{\mathcal{X}}$ is highly complex, we can resort to the numerical methods such as grid-based numerical integration \cite{zhong2023efficient}.

Note that \eqref{SOS1} is a bilinear semi-definite problem which is hard to solve \cite{henrion2014convex}.
We can empirically provide the initial value of $\lambda$ in advance to avoid the bilinear problem, which ultimately aligns with the framework of convex optimization, making the problem computationally feasible. 
Then, we can efficiently solve the program \eqref{SOS1} in polynomial time using
the interior point method, leveraging the SOS module YALMIP \cite{lofberg2004yalmip} and the Mosek optimization toolbox \cite{aps2019mosek}.

\subsection{ROS around Equilibrium Point}
\label{3-2}

In this subsection, we consider a special case where the synchronization manifold $\mathcal{S}$ becomes an equilibrium point $\mathbf{S}=[0^{\top}_{n},0^{\top}_{n},\ldots,0^{\top}_{n}]^{\top} \in \mathbb{R}^{nN}$, as in \cite{zhu2019estimating}. Then, network \eqref{network} is reformulated in a compact form as
\begin{align} \label{network2}
\dot{\mathbf{X}}=\mathbf{F}(\mathbf{X})-c(L \otimes I_{n}) \mathbf{G}(\mathbf{X}),
\end{align}
where $\mathbf{F}(\mathbf{X})=[f^{\top}(\mathbf{x}_{1}),f^{\top}(\mathbf{x}_{2}), \ldots,f^{\top}(\mathbf{x}_{N})]^{\top} \in \mathbb{R}^{nN}$ and $\mathbf{G}(\mathbf{X})=[g^{\top}(\mathbf{x}_{1}),g^{\top}(\mathbf{x}_{2}), \ldots,g^{\top}(\mathbf{x}_{N})]^{\top} \in \mathbb{R}^{nN}$.

{\bf Remark 4}. Once synchronization is achieved, the coupling term $c(L \otimes I_{n}) \mathbf{G}(\mathbf{X})$ vanishes, and the evolution of the synchronized trajectories is governed solely by the intrinsic node dynamics $f(\cdot)$. 
If the node dynamics possesses a stable equilibrium point, the synchronized trajectories will converge to it, i.e., the synchronization manifold becomes an equilibrium point (see Example 2).
Conversely, if the node dynamics exhibits a limit cycle or other oscillatory attractor, the synchronization manifold will reflect that periodic (or possibly chaotic) behavior (see Example 1).

Analogous to Section \ref{3-1}, we can naturally obtain the following theorem and optimization program.

{\bf Theorem 2}. If there exists a continuously differentiable function $V(\mathbf{X}): \overline{\mathcal{X}} \rightarrow \mathbb{R}$ that satisfies
\begin{align}
\mathcal{L}_{V}(\mathbf{X})-\lambda V(\mathbf{X}) \geq 0, \qquad
&\ \forall \mathbf{X} \in \overline{\mathcal{X} \setminus \mathcal{X}_{T}}, \label{EGBF3} \\
V(\mathbf{X}) \leq 0, \qquad \qquad \qquad \ \
&\ \forall \mathbf{X} \in  \partial \mathcal{X} \label{EGBF4}
\end{align}
with $\lambda > 0$, then $\mathcal{R}=\{\mathbf{X} \in \mathcal{X} \mid V(\mathbf{X})>0 \}$ is a subset of ROS around equilibrium point for network \eqref{network2}.

The proof follows the lines of the proof of Theorem 1 and is omitted.

The SOS program to be solved for $V(\mathbf{X})$ is as follows:
\begin{equation} \label{SOS2}
\begin{split}
&\max_{\boldsymbol{c},\lambda} \ \boldsymbol{c}^{\top}  \boldsymbol{\rho} \\
&\ \text{s.t.}  \ \mathcal{L}_{V}(\mathbf{X})-\lambda V(\mathbf{X})
- p_{1}(\mathbf{X})h(\mathbf{X}) \\
&\quad \quad \quad \quad \quad \quad \quad \quad \quad + p_{2}(\mathbf{X})l(\mathbf{X}) \in \sum[\mathbf{X}],    \\
& \quad \  -V(\mathbf{X})+q(\mathbf{X})h(\mathbf{X}) \in \sum[\mathbf{X}], \\
& \quad \ \ p_{1}(\mathbf{X}) \in \sum[\mathbf{X}],
\ p_{2}(\mathbf{X}) \in \sum[\mathbf{X}], \ q(\mathbf{X}) \in \mathbb{R}[\mathbf{X}],
\end{split}
\end{equation}
where $\boldsymbol{c}^{\top} \boldsymbol{\rho} = \int_{\overline{\mathcal{X}}} V(\mathbf{X})d(\mathbf{X})$, $\boldsymbol{\rho}$ is a constant vector computed by integrating the monomials in $V(\mathbf{X}) \in \mathbb{R}[\mathbf{X}]$ over $\overline{\mathcal{X}}$, and $\boldsymbol{c}$ is a vector composed of unknown coefficients in $V(\mathbf{X})$; and $l(\mathbf{X})=\varepsilon - ||\mathbf{X}||$.

{\bf Remark 5}. A typical result for estimating the ROS for network \eqref{network2} takes the form of a spheroid \cite{zhu2019estimating}, where the center corresponds to the equilibrium point $\mathbf{S}$, and the radius $r_{0}$ represents the size of the estimated ROS (i.e., $||\mathbf{X}||<r_{0}$).
It is well known that the shape of the ROS is generally complex, and the method proposed in \cite{zhu2019estimating} is conservative, as it estimates the ROS using a simple spheroid.
In contrast, our approach generates a ROS $\mathcal{R}$ with a more general shape. Furthermore, the size of $\mathcal{R}$ is larger compared to the spheroid approximation (see Example 2 for more details).

\vspace{0.5em}
\section{Numerical Examples}
\label{4}

In this section, we present two numerical examples to illustrate the validity of our proposed method.

{\bf Example 1}. Consider a networked system \eqref{network} with three nodes, each of which is described as the following Van der Pol oscillator \cite{lee2018heterogeneous}:
\begin{align}
\left\{\begin{array}{lr}
\dot{x}_{i}=y_{i} \\
\dot{y}_{i}=\hat{\mu}\hat{\omega}(1-\hat{r}x_{i}^{2})y_{i}-\hat{\omega}^2x_{i}^{2}
\end{array}\right. \label{vandepol}
\end{align}
where $\hat{r}=1$, $\hat{\mu}=0.5$ and $\hat{\omega}=0.9$.
Let $L=[4,-2,-2;-1,2,1;-3,0,3]$ be the Laplacian matrix of digraph $\mathscr{G}$,
and choose the coupling strength $c=0.1$ and the linear coupling function $g(\mathbf{x}_{i})=[x_{i},y_{i}]^{\top}$ with $i=1, 2, 3$.

\begin{figure}[!ht]
\centering
{\includegraphics[width=3in]{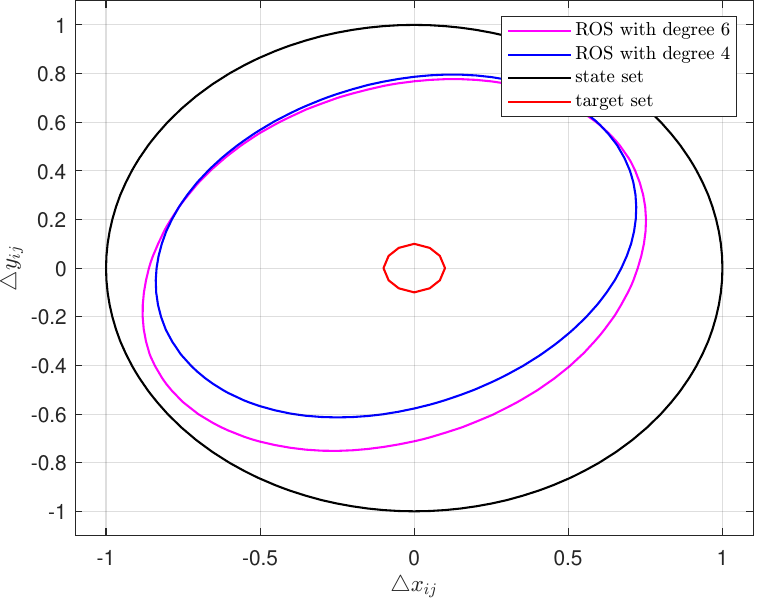}}
\caption{Black and red curves -- $\partial \mathcal{X}$ and $\partial \mathcal{X}_{T}$; blue and purple curves -- $\partial \mathcal{R}$ in Example 1.}
\label{EX1}
\end{figure}

\begin{figure}[t]
\centering
\subfigure
{		\includegraphics[width=2.5in]{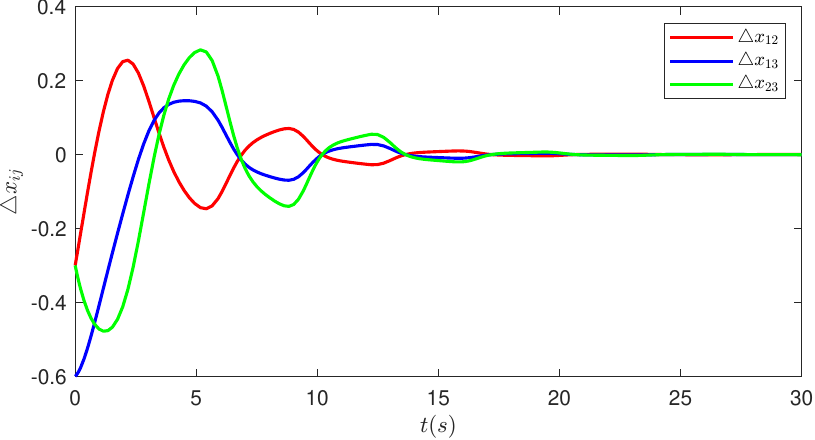}}
\subfigure
{	    \includegraphics[width=2.5in]{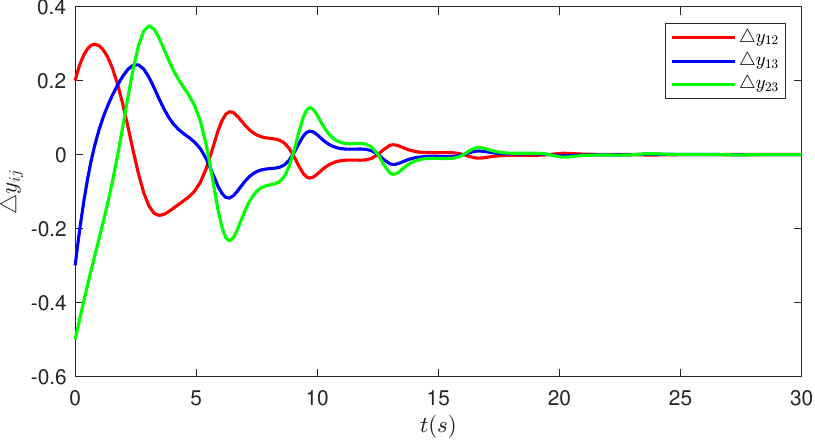}}
\caption{The error trajectories of network \eqref{vandepol}.}
\label{EX1-2}
\end{figure}

It follows from \cite{wang2021decomposition} that the synchronization manifold $\mathcal{S}$ of network \eqref{vandepol} is stable. Given the state set $\mathcal{X}=\{(\triangle x_{ij},\triangle y_{ij})^{\top} | \triangle x_{ij}^{2}+\triangle y_{ij}^{2}<1\}$
and the target set $\mathcal{X}_{T}=\{(\triangle x_{ij},\triangle y_{ij})^{\top} | \triangle x_{ij}^{2}+\triangle y_{ij}^{2}<0.01\}$ for $i,j = 1, 2, 3$ and $i \neq j$,  the ROS $\mathcal{R}$ is computed for degrees 4 and 6 of the EGBF $V$ by solving the SOS program \eqref{SOS1} with $\lambda=0.1$. The computed results are illustrated in Fig. \ref{EX1}.
Moreover, for the initial values $[0.2,0,0.5,-0.2,0.8,0.3]^{\top} \in \mathcal{R}$, Fig. \ref{EX1-2} shows that the error trajectories of network \eqref{vandepol} converge to zero, which validates the theoretical results.

{\bf Example 2}. Consider a coupled model of Moore-Greitzer jet engines in the no-stall mode \cite{zhang2023stability}. The networked system \eqref{network} consists of four jet engines, each described by
\begin{align}
\left\{\begin{array}{lr}
\dot{x}_{i}=-0.5{x}_{i}^{3}-1.5{x}_{i}^{2}-y_{i} \\
\dot{y}_{i}=3{x}_{i}-{y}_{i}
\end{array}\right. \label{jet-engines}
\end{align}
where $x_{i}$ relates to the mass flow and $y_{i}$ relates to the pressure rise.
Let $L=[2,-1,-1,0;-1,1,0,0;0,0,1,-1;$ $-1,0,0,1]$ be the Laplacian matrix of digraph $\mathscr{G}$, and choose the coupling strength $c=1$ and the nonlinear coupling function $g(\mathbf{x}_{i})=[x_{i},y_{i}+0.5y_{i}^{3}]^{\top}$ with $i=1, 2, 3, 4$.

\begin{figure}[t]
\centering
{\includegraphics[width=3in]{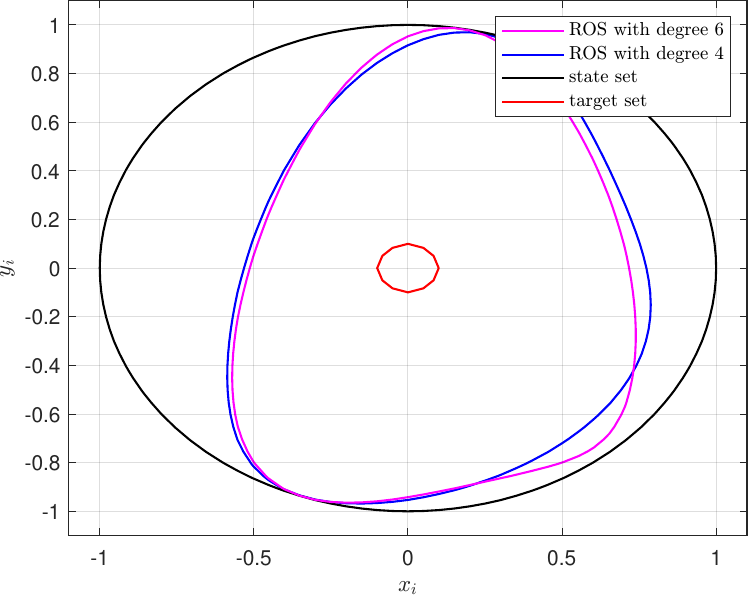}}
\caption{Black and red curves -- $\partial \mathcal{X}$ and $\partial \mathcal{X}_{T}$; blue and purple curves -- $\partial \mathcal{R}$ in Example 2.}
\label{EX2}
\end{figure}

\begin{figure}[!ht]
\centering
\subfigure
{		\includegraphics[width=2.5in]{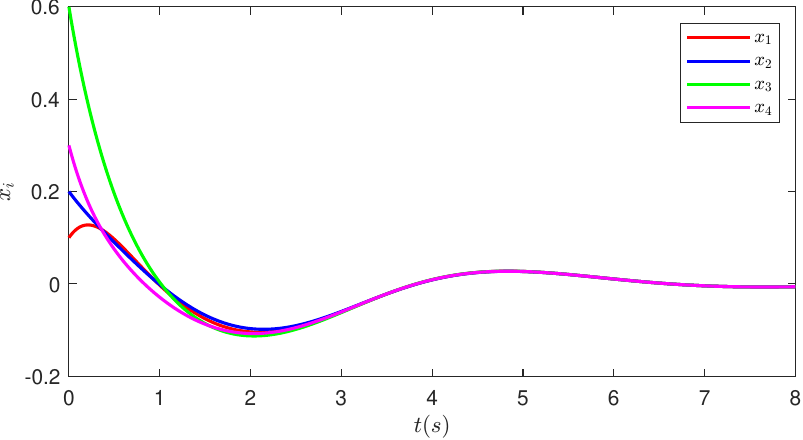}}
\subfigure
{	    \includegraphics[width=2.5in]{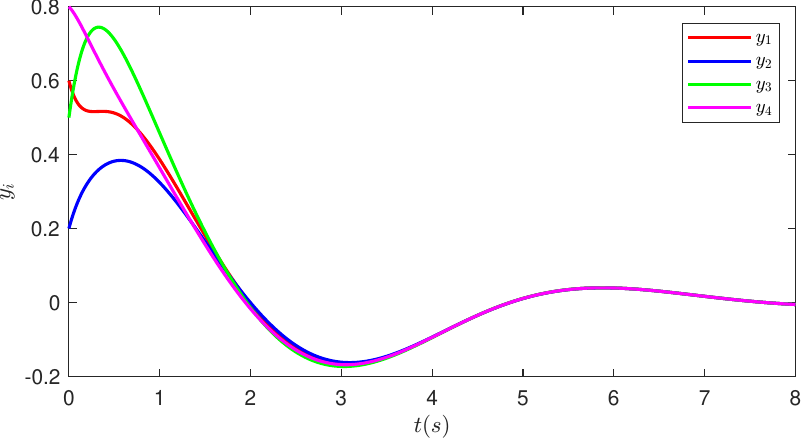}}
\caption{The state trajectories of network \eqref{jet-engines}.}
\label{EX2-2}
\end{figure}

According to \cite{zhang2023stability}, the synchronization manifold of network \eqref{jet-engines} is a stable equilibrium point $\mathbf{S}=[0,0,0,0,0,$ $0,0,0]^{\top}$.
Given the state set $\mathcal{X}=\{(x_{i}, y_{i})^{\top} | x_{i}^{2}+ y_{i}^{2} < 1\}$
and the target set $\mathcal{X}_{T}=\{(x_{i}, y_{i})^{\top} | x_{i}^{2} + y_{i}^{2}<0.01\}$ for $i= 1, 2, 3, 4$, 
the ROS $\mathcal{R}$ is computed for degrees 4 and 6 of the EGBF $V$ by solving the SOS program \eqref{SOS2} with $\lambda=0.01$ (see Fig. \ref{EX2}). 
Apparently, our results outperform the estimated ROS ($||\mathbf{X}|| < 0.0238$) in \cite{zhu2019estimating}, yielding a larger ROS with a more general shape.
Besides, for the initial values $[0.1,0.6,0.2,0.2,0.6,0.5,0.3,0.8]^{\top} \in \mathcal{R}$, Fig. \ref{EX2-2} shows that the state trajectories of network \eqref{jet-engines} converge to zero, which illustrates the effectiveness of our approach.

\vspace{1em}
\section{Conclusion}
\label{5}

This article has investigated the ROS estimation of complex networks based on the SOS decomposition for multivariate polynomials.
The core of our method was the introduction of EGBF to construct the corresponding sufficient conditions for generating the ROS along the synchronization manifold.
The ROS was then efficiently computed by solving a convex SOS programming problem.
Moreover, the proposed method has also been applied to estimate the ROS around the equilibrium point, allowing it to take more general shapes.
Finally, the numerical examples have been provided to showcase the effectiveness of our method.

In the future, we would extend the proposed SOS programming approach to address the ROS problems for heterogeneous networks with nonidentical dynamics, building on the foundations laid by earlier works \cite{wang2021decomposition,zhang2023stability,zhang2024consensus}.

\addtolength{\textheight}{-8cm}


\vspace{1em}
\bibliographystyle{ieeetr}
\bibliography{ref_region}
\end{document}